\documentclass[final]{siamltex}

\usepackage{graphicx, amssymb, amsmath}
\usepackage{color}
\usepackage{pst-text}
\usepackage{soul} 
\usepackage{hyperref}

\graphicspath{{./}{figs/}}

\newcommand{\R}{\mathbb R}

\newcommand{\approxT}{\widehat T}
\newcommand{\approxH}{\widehat H}
\newcommand{\ACC}{\textit{ACC}}
\newcommand{\CONS}{\textit{CONS}}
\newcommand{\FAR}{\textit{FAR}}
\newcommand{\ACCsp}{\textit{ACC }}
\newcommand{\CONSsp}{\textit{CONS }}
\newcommand{\FARsp}{\textit{FAR }}
\newcommand{\be}{\begin{equation}}
\newcommand{\ee}{\end{equation}}

\newtheorem{dfn}{Definition}[section]

\newtheorem{rmk}{Remark}[section]

\title{{Can local single--pass methods solve any stationary Hamilton--Jacobi--Bellman equation?}\thanks{The authors wish to acknowledge the support obtained by the following grants: AFOSR Grant no.\ FA9550-10-1-0029, 
ITN - Marie Curie Grant no.\ 264735-SADCO 
and SAPIENZA 2009 {``Analisi ed approssimazione di modelli differenziali nonlineari in fluidodinamica e scienza dei materiali''}.}
}

\author{
Simone Cacace\footnote{Dipartimento di Matematica, SAPIENZA -- Universit\`a di Roma, P.le Aldo Moro 2, 00185 Rome, Italy. 
\texttt{cacace@mat.uniroma1.it}}, 
Emiliano Cristiani\footnote{Istituto per le Applicazioni del Calcolo, Consiglio Nazionale delle Ricerche, Rome, Italy. 
\texttt{e.cristiani@iac.cnr.it}},
Maurizio Falcone\footnote{(corresponding author) Dipartimento di Matematica, SAPIENZA -- Universit\`a di Roma, P.le Aldo Moro 2, 00185 Rome, Italy. 
\texttt{falcone@mat.uniroma1.it}} 
}

\begin{document}
\maketitle

\begin{abstract}
The use of \textit{local single-pass} methods (like, e.g., the Fast Marching method) has become popular in the solution of some Hamilton-Jacobi equations. 
The prototype of these equations is the \textit{eikonal equation}, for which the methods can be applied saving CPU time and possibly memory allocation. 
Then, some natural questions arise: can local single-pass methods solve any Hamilton-Jacobi equation? If not, where the limit should be set?\\ 
This paper tries to answer these questions. 
In order to give a complete picture, we present an overview of some fast methods available in literature and we briefly analyze their main features. 
We also introduce some numerical tools and provide several numerical tests which are intended to exhibit the limitations of the methods. 
We show that the construction of a local single-pass method for general Hamilton-Jacobi equations is very hard, if not impossible. Nevertheless, some special classes of problems can be actually solved, making local single-pass methods very useful from the practical point of view.
\end{abstract}



\pagestyle{myheadings}
\thispagestyle{plain}
\markboth{S. CACACE, E. CRISTIANI, M. FALCONE}
{CAN LSP METHODS SOLVE ANY STATIONARY HJB EQUATION?}


\section{Introduction}\label{sec:intro}
The study of Hamilton-Jacobi (HJ) equations arises in several applied contexts, including classical mechanics, front propagation, control problems and differential games, and it has a  great impact 
in many areas, such as  robotics, aeronautics, electrical and aerospace engineering. 
In particular, for control/game problems, an approximation of the value function allows for the synthesis of optimal control laws in feedback form, 
and then for the computation of optimal trajectories. The value function for a control problem (resp., differential game) can be characterized as the unique viscosity solution 
of the corresponding Hamilton-Jacobi-Bellman (HJB) equation (resp., Hamilton-Jacobi-Isaacs (HJI) equation), and it is obtained by passing to the limit in the well known Bellman's Dynamic Programming (DP) principle. 
The DP approach can be rather expensive from the computational point of view, but in some situations it gives a real advantage when compared to methods based on the Pontryagin's Maximum Principle, because the latter approach allows one to compute only open-loop controls and locally-optimal trajectories. 
Moreover, weak solutions to HJ equations are nowadays well understood in the framework of viscosity solutions, which offers the correct notion of solution for many applied problems.

The above remarks  have motivated  the research of efficient and accurate numerical methods. Indeed, an  increasing number of techniques have been proposed for the approximation of viscosity solutions. 
They range from Finite Difference to Finite Volume, from Discontinuous Galerkin to semi-Lagrangian schemes. In any case, for optimal control problems and games, the DP approach suffers from the so-called ``curse of dimensionality'' limitation, i.e.\ the size of  nonlinear systems to solve for high dimensional problems becomes huge and this makes the numerical solution unfeasible, in terms of both memory allocation and CPU time.  The curse of dimensionality can be sometimes overcome by exploiting the peculiarities of the problem, if any (e.g., symmetry, periodicity, linearity), or by adopting a linearization based on the so-called ``max-plus algebra'' approach, which unfortunately presents other types of  constraints, see, e.g., the book by McEneaney \cite{McE06}. 
It is rather clear that the DP approach needs a big effort in the construction of numerical approximation schemes for two different reasons. The first, which is valid even in low dimension, is due to the low regularity of viscosity solutions which are typically only Lipschitz continuous or even discontinuous, as in the case of constrained control problems and pursuit-evasion games. 
The second reason is related to the above mentioned curse of dimensionality which pushes towards methods with low memory allocation and, possibly, the definition of some rule to reduce the number of elementary operations and the CPU time. 

Another motivation for efficient numerical methods is the approximation of front propagation problems via the level-set method. The motivation there is to reduce or eliminate the extra dimension which is added by the level-set method and obtain a fast and reliable algorithm.
Starting from the 1980s, many efforts have been made to improve the efficiency of these numerical methods, a crucial step for the solution of real-world problems.

\medskip
In this paper we deal with numerical methods for solving first-order nonlinear convex stationary HJB equations. In particular, we focus on the applicability of 
Fast Marching Method (FMM), introduced in the pioneering works by Tsitsiklis \cite{T95}, Sethian \cite{S96}, Helmsen et al. \cite{HPCD96}, and its generalizations, see, e.g., \cite{AM12, CFFM08, CV12, C01, C09, CF08, CF07, JW08, K01, PS05, SV03}. 
We analyze features and limitations of this kind of algorithms, aiming at understanding if it is possible to construct local single-pass methods  (see definitions \ref{def1} and \ref{def2} below) solving every HJB equation. Then, we discuss whether or not the research on this topic should look for new future directions still based on the local single-pass idea and/or switch to other acceleration methods, such as Fast Sweeping Methods (FSM), see, e.g., \cite{BD99, KOQ04, TCOZ04, Z05, Z00} (see also the pioneering work \cite[p.168]{K77}, where a sketch of the method is given, and \cite{D80}, where a similar method in a discrete setting is proposed).

It is well known that FMM is an efficient numerical technique for solving the eikonal equation. This explains why we decide to use as a guide line the following equations, which generalize the eikonal equation and are associated to some minimum time problems with target:
\begin{eqnarray}
 & &\sup\limits_{a\in B(0,1)}\left\{-a\cdot\nabla T(x)\right\} = 1
 \qquad \qquad \phantom{'} \textrm{(homogeneous eikonal)}
 \label{eiko_omog}
 \\ 
 & &\sup\limits_{a\in B(0,1)}\left\{-c_1(x)a\cdot\nabla T(x)\right\} = 1
 \qquad \textrm{(nonhomogeneous eikonal)}
 \label{eiko_nonomog}
 \\ 
 & &\sup\limits_{a\in B(0,1)}\left\{-c_2(a)a\cdot\nabla T(x)\right\} = 1
 \qquad \textrm{(homogeneous anisotropic eikonal)}
 \label{eiko_aniso_omog}
 \\
 & &\sup\limits_{a\in B(0,1)}\left\{-c_3(x,a)a\cdot\nabla T(x)\right\} = 1
 \quad  \textrm{(nonhomogeneous anisotropic eikonal)}
 \label{eiko_aniso_nonomog}
 \\
 & &\sup\limits_{a\in B(0,1)}\left\{-f(x,a)\cdot\nabla T(x)\right\} = 1
 \qquad \! \textrm{(minimum time HJB)}
 \label{sect1:HJB}
\end{eqnarray}
where $x\in\R^d\backslash\mathcal T$, $\mathcal T$ is a closed nonempty target set in $\R^d$, $c_1$, $c_2$, $c_3$ are given strictly positive and Lipschitz continuous scalar functions, $f$ is a given vector-valued Lipschitz continuous function, and $B(0,1)$ is the unit ball in $\R^d$, representing the set of the admissible controls. 
To simplify the presentation we will always consider the homogeneous Dirichlet condition $T=0$ on $\mathcal T$ but also other boundary conditions can be applied, 
provided some compatibility conditions between the vectorfield $f$ and $\partial\mathcal T$ hold true. 
Let us also note, for the readers not familiar with control applications, that equations (\ref{eiko_omog}) and (\ref{eiko_nonomog}) can be written in a more standard form as
$$
|\nabla T(x)|=1 \qquad \text{ and } \qquad c_1(x)|\nabla T(x)|=1,
$$ 
respectively.
Moreover, the above relation shows the equivalence between the front propagation problem described by the level set method and the minimum time problem, as one can find in \cite{F94}.
To simplify the notations, we restrict the discussion to the case $d=2$, but the results of the paper are valid in any dimension.

It is interesting to note that the single-pass idea has been also applied to nonconvex hamiltonians, e.g.\ the HJI equation corresponding to pursuit-evasion games \cite{CCF11, CF06}. In pursuit-evasion games the structure of the solution is similar to the minimum time problem because the characteristic information propagates from a given target to the rest of the space. Clearly for games  the structure of optimal trajectories (which coincide with the characteristics of the problem) is much more complicated due to the presence of two independent players. However, a complete theory for Fast-Marching-like methods in this framework is still missing and, as we will see later,  it will be hard to develop it since  these methods can fail even in the convex case  (see the examples in Section \ref{sec:numtests}).

As mentioned above, in the last decades many numerical schemes and algorithms were proposed to solve the above equations. 
Some of these schemes are listed in the next section, together with their main properties.  
As it is well known, one important feature held by Fast-Marching-like methods is that the solution to the HJ equation is computed in a finite number of steps.  
More precisely, these methods are \emph{single-pass}, in the sense of the following definition.

\begin{dfn}[Single-pass algorithm]\label{def1}
An algorithm is said to be \emph{single-pass} if each mesh point is re-computed at most $r$ times, where $r$ depends only on the equation and the mesh structure, 
not on the number of mesh points.
\end{dfn}

Single-pass algorithms usually divide the numerical grid in, at least, three subsets: \emph{Accepted} (\ACC) region, \emph{Considered} (\CONS) region and \emph{Far} (\FAR) region. 
Nodes in \ACCsp are definitively computed, nodes in \CONSsp are computed but their values are not yet final, and nodes in \FARsp are not yet computed.

We also introduce the following definition.

\begin{dfn}[Local single-pass algorithm]\label{def2}
A single-pass algorithm is said to be \emph{local} if the computation at any mesh point involves only the values of first neighboring nodes, 
the region \CONSsp is one-cell thick  and no information coming from \FARsp region is used.
\end{dfn}

\medskip

The paper is organized as follows: 
in Section 2 we summarize some of the existing methods to solve HJB equations and introduce two semi-Lagrangian numerical schemes. In Section 3, we present new numerical tools which will be useful to investigate the applicability of local single-pass methods. In Section 4, which is the core of the paper, 
we discuss the main features and limitations of the methods presented in Section 2 and 3, and we address the problem of extending local single-pass methods to general HJB equations. 
Finally, in Section 5  we present several experiments and numerical tests, in order to compare the two schemes described in Section 2 and to confirm the scenario depicted in Section 4.


\section{Background and general approximation schemes}
\label{sec:background}
Fast methods for HJB equations are usually designed to work with different local schemes, including Finite Difference and semi-Lagrangian (SL) schemes. Several results show that, in many cases on structured grids and at a reasonable cost, SL schemes provide better accuracy than other schemes (see, e.g., \cite{CF07, F97}), due to their ability to follow directions which are oblique with respect to the coordinate axes.
In this section we recall, for readers convenience, two SL schemes for HJB equations, which will be compared in Section \ref{sec:numtests}. 
Then, we list and discuss some of the iterative and Fast-Marching-like methods available in literature.
\subsection{Two SL schemes} \label{numschemes}
Let us  introduce a structured grid $G$ and denote its nodes by $x_i$, $i=1,\ldots,N$. 
The space step is assumed to be uniform and equal to $\Delta x>0$. HJ equations can be discretized by means of the discrete version of the DP principle. 
In this way the relationship with the optimal control framework is never lost. Standard arguments \cite{BCDbook} lead to the following discrete version of the HJB equation (\ref{sect1:HJB}): 

\begin{equation}\label{SLscheme}
T(x_i)\approx \approxT(x_i) = 
\min\limits_{a\in B(0,1)}\left\{ \approxT(\tilde x_{i,a})+ \frac{|x_i-\tilde x_{i,a}|}{|f(x_i,a)|}\right\}\,,\qquad x_i\in G
\end{equation}
where $\tilde x_{i,a}$ is a \textit{non-mesh} point, obtained by integrating, until a certain final time $\tau$, the ordinary differential equation
\be \label{ODE}
\left\{
\begin{array}{ll}
\dot y(t)=f(y,a), & t\in[0,\tau] \\
y(0)=x_i & 
\end{array} 
\right .
\ee
and then setting $\tilde x_{i,a}=y(\tau)$. To make the scheme fully discrete, the set of admissible controls $B(0,1)$ is discretized in $N_c$ points and we denote by $a^*$ the optimal control 
achieving the minimum in (\ref{SLscheme}).
Note that we can get different versions of the SL scheme (\ref{SLscheme})
varying $\tau$, the method used to solve (\ref{ODE}), and the interpolation method used to compute  $\approxT (\tilde x_{i,a})$. 
Moreover, we remark that, in any single-pass method, the computation of $\approxT (x_i)$ cannot involve the value $\approxT (x_i)$ itself, because this self-dependency would make the method iterative.
\subsubsection*{A two-point SL scheme}
\begin{figure}[h!]
\begin{center}
\begin{tabular}{ccc}
\includegraphics[width=.35\textwidth]{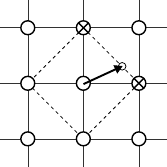} & \qquad \qquad &
\includegraphics[width=.35\textwidth]{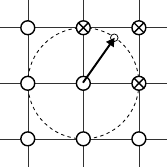} \\
(a) & \qquad\qquad & (b)
\end{tabular}
\put(-253.5,25){\small $f(x_i,a)$}
\put(-257,0){\small $x_i$}
\put(-198,0){\small $x_{i,1}$}
\put(-262,64){\small $x_{i,2}$}
\put(-88,34){\small $f(x_i,a)$}
\put(-83,0){\small $x_i$}
\put(-24,0){\small $x_{i,1}$}
\put(-88,64){\small $x_{i,3}$}
\put(-24,64){\small $x_{i,2}$}
\caption{(a) 2-points SL scheme. (b) 3-points SL scheme}
\label{fig:SL2p3p}
\end{center}
\end{figure}
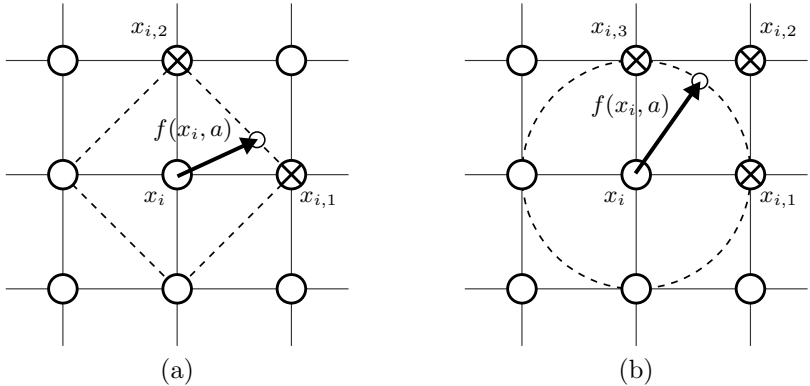

This scheme is used, for example, in  \cite{SV03} and \cite{T95}. Equation (\ref{ODE}) is solved by an explicit forward Euler scheme until the solution intercepts the line connecting 
two neighbouring points $x_{i,1}$ and $x_{i,2}$ (see Fig.\ \ref{fig:SL2p3p}a).
The value $\approxT (x_i)$ is computed by a one-dimensional linear interpolation of the values $\approxT (x_{i,1})$ and $\approxT (x_{i,2})$ with weights $\lambda_{i,1}$ and $\lambda_{i,2}$ respectively 
($\lambda_{i,1}+\lambda_{i,2}=1$).
%

\subsubsection*{A three-point SL scheme}
This scheme is used, for example, in \cite{CF07}. Equation (\ref{ODE}) is solved by an explicit forward Euler scheme until the solution is at distance $\Delta x$ from $x_i$, 
where it falls inside the triangle of vertices $x_{i,1}$, $x_{i,2}$, and $x_{i,3}$ (see Fig.\ \ref{fig:SL2p3p}b). 
The value $\approxT (x_i)$ is computed by a two-dimensional linear interpolation of the values $\approxT (x_{i,1})$, $\approxT (x_{i,2})$ and $\approxT (x_{i,3})$ 
with weights $\lambda_{i,1}$, $\lambda_{i,2}$ and $\lambda_{i,3}$ respectively ($\lambda_{i,1}+\lambda_{i,2}+\lambda_{i,3}=1$).\\

\begin{rmk}\label{rem:includeACC}
It is important to note that, when the algorithm employs the single-pass technique (with \ACC-\CONS-\FARsp sets), the nodes in \CONSsp can be either included in or excluded from the computation of $\approxT(x_{i})$. 
Indeed, we can decide to force the scheme to employ only nodes in \ACC, temporarily assuming that nodes in \CONSsp have very large values, so that they are automatically rejected by the minimum search
 in \eqref{SLscheme}. Otherwise, we can employ nodes in \ACCsp and \CONS, although the values at nodes in \CONSsp are not in general correct since they can still vary in the following iterations. 
\end{rmk}
\subsection{Some algorithms for HJB equations}\label{algorithms}
Here we list and briefly describe some iterative and single-pass methods for solving HJB equations.

\medskip
\paragraph{Iterative Method (ITM)} 

This method naturally exploits the fixed point form of the discrete DP principle. It has been applied in \cite{K77} to stochastic control problems, where the deterministic case can be obtained vanishing the coefficient in front of the diffusion term. Later, a similar operator has been used in \cite{GR85} to construct an approximation of the value function in a deterministic control problem. In the framework of viscosity solutions, the fully discrete DP scheme was first studied in \cite{F87}. The interested reader can also find a detailed description of the algorithm and some acceleration methods in \cite{F97}.
A finite difference version for a specific application to the Shape-from-Shading problem can be found in \cite{RT92}.

Starting from some initial guess $\approxT^{(0)}$ defined on the whole grid (compatible with the Dirichlet conditions imposed on $\mathcal T$) the nodes are visited in some unique and predefined order. At each visit, the numerical scheme is applied and a new value for the node is computed. 
This leads to a fixed-point algorithm of the form
$$
\approxT^{(n)}=\approxH(\approxT^{(n-1)})\,,\qquad n=1,2,3,\ldots\,,
$$
where $\approxH$ denotes a discrete Hamiltonian associated to the corresponding HJ equation. 
Gauss-Seidel-type or Jacobi-type iterations are possible. For a practical implementation, a criterion of the form
\be\label{stopcriterion}
\max_{x_i\in G}|\approxT^{(n)}(x_i)-\approxT^{(n-1)}(x_i)|<tol
\ee
is needed in order to stop the computation at a desired precision \emph{tol}. Clearly this method is not single-pass, since the number of iterations needed to reach convergence depends both on the grid size $\Delta x$ and the dynamics underlying the equation. 
ITM was proved to be convergent, provided a suitable numerical scheme is employed.

\medskip
\paragraph{Fast Sweeping Method (FSM) \cite{BD99,KOQ04,TCOZ04,Z05,Z00}} This method is similar to ITM, but the grid is visited in a multiple-direction predefined order. 
Usually, a rectangular grid is iteratively swept along four directions: $N\rightarrow S$, $E\rightarrow W$, $S\rightarrow N$ and $W\rightarrow E$, where $N,S,E,W$ stand for North, South, East and West respectively. 
This method has been shown to be much faster than ITM, but, as ITM, it not single-pass.
A well known exception is given by the eikonal equation, for which it is proved that only $1$ sweep (i.e. four visits of the whole grid) is enough to reach convergence (see \cite{Z05} for details). 
FSM computes the same solution of ITM, provided the same scheme and the same stopping rule are employed. 

\medskip
\paragraph{Fast Marching Method (FMM) \cite{HPCD96,S96,T95}} This method has been introduced as a fast solver for the eikonal equation. It differs from the previous ones, since the nodes are visited in a solution-dependent order, producing a single-pass method: 
the algorithm itself finds a correct order for processing the grid nodes. 
The order which is determined satisfies the \textit{causality} principle, i.e.\ the computation of a node is declared completed only if its value cannot be affected by the future computation. As recalled in Section \ref{sec:intro}, at each step the grid is divided in three regions: 
\ACC, where computation is definitively done, \CONS, where computation is going on and \FAR, where computation is not done yet. 
Then, the node in \CONSsp with the minimal value  enters \ACC, its first neighbours enter \CONSsp (if not yet in) and are (re)computed.\\
Following \cite{SV03}, we remark that this \emph{minimum-value rule} corresponds to compute the value function $T$ step by step in the ascending order (i.e., from the simplex containing $-\nabla T$). 
It follows that \CONSsp expands under the gradient flow of the solution itself, which is exactly equivalent to say that \CONSsp is, at each step, an approximation of a level set of the value function. 
In the case of isotropic eikonal equation (\ref{eiko_nonomog}), the gradient of the solution coincides with the characteristic field, hence FMM computes the correct solution. 
Moreover, FMM still works for problems with mild anisotropy, where gradient lines and characteristics define small angles and lie, at each point, in the same simplex of the underlying grid. 
On the other hand, when a strong anisotropy comes into play, as for a general  anisotropic eikonal equation (\ref{eiko_aniso_omog}), 
FMM fails and there is no way to compute the viscosity solution following its level sets.   
Finally, we remark that FMM is also a local method, since each node is computed by means of first neighbors nodes only and \CONSsp is one-cell thick.  Moreover, FMM computes the same solution of ITM, provided the same scheme is employed. 

\medskip
\paragraph{Characteristic Fast Marching Method (CFMM) \cite{CF08}}
This method is inspired by FMM, it is local and single-pass and can be used to solve some eikonal equations. 
It replaces the search for the minimum value in \CONSsp with the search of the node where the characteristic line passes with maximal speed. 
The acceptance rule is also modified: a node $x_i$ in \CONSsp enters \ACCsp if the point $x_i+f(x_i,a^*)$ falls in \ACC.
As the Group Marching Method \cite{K01}, more than one node can enter \ACCsp at the same time, making the method in general faster than FMM. 
Note that CFMM does not always work if the solution of the equation is not differentiable.

\medskip
\paragraph{Ordered Upwind Method (OUM) \cite{SV03}} This method is inspired by FMM, but it is able to solve more general equations than the eikonal one, 
including nonhomogeneous anisotropic eikonal equations (\ref{eiko_aniso_nonomog}). 
This can be obtained by enlarging the stencil of the scheme, so that a value at a node $x_i$ can be computed by using values at some nodes $x_j$ that are 
far from the node $x_i$. This makes the method nonlocal. The maximal allowed distance $|x_i-x_j|$ depends on the degree of anisotropy of the equation. 
OUM is a single-pass method which computes the same solution of ITM (employing the same numerical scheme) only in the limit $\Delta x\rightarrow 0$.\\ 
A generalization of OUM has been recently proposed in \cite{AM12} to solve static convex HJ equations on highly nonuniform grids. 
The new method MAOUM (Monotone Acceptance OUM) computes the solution in a Fast Marching fashion, but employs large stencils (even larger than OUM) that are pre-computed for each grid node. 
This makes MAOUM two-pass and nonlocal. 

\medskip
\paragraph{Buffered Fast Marching Method (BFMM) \cite{C09}} This method is inspired by FMM and can be used to solve any HJ equations modelling monotone front propagation. 
Although only first neighbours are involved in the computation, BFMM cannot be considered a local method, since \CONSsp can increase its thickness. 
More precisely, the \CONSsp region is extended by a \emph{Buffer} region, whose size depends on the dynamics of the equation and, in the worst case, it can cover the whole grid, thus making BFMM comparable to ITM.
BFMM is not single-pass and computes the same solution of ITM, provided the same scheme is employed.

\medskip
\paragraph{Progressive Fast Marching Method (PFMM) \cite{CCF11}} This method can be considered a localization of BFMM. It is indeed a local method, but not single-pass. Some experimental results have shown that 
it can solve quite 
general problems, including Pursuit-Evasion games with state and control constraints. PFMM has been introduced for theoretical purposes only, since it is very slow (slower than ITM) and then not usable in practice. 
It proposes a completely new rule for accepting nodes in \CONS: in the \FARsp region, next to the \CONSsp region, a layer of ``tempting'' values is placed and \emph{progressively} increased, 
acting as an external boundary condition. For each tempting value, the solution is re-computed in \CONS, recording the corresponding variations. 
The first node in \CONSsp which is not affected by this external layer enters \ACC. The ``tempting'' values can be considered as a guess on the outcome of the future computation
and the new rule of acceptance allows one to find the node in \CONSsp that cannot be affected by it.

\section{New tools and verification methods}\label{sec:newtools}

In this section we consider four additional FM-like methods. 
The first two are acceleration techniques which are expected to provide the same solution of ITM whenever they are applicable.
The last two, labelled as \emph{dumb}, are not new methods for solving HJ equations, rather they are \textit{verification tools}. They will be used to analyze features and limitations of the methods already presented, 
aiming at giving a comprehensive classification of the equations that can be solved by local single-pass algorithms. 
Our ultimate goal is to discuss the possibility that local single-pass methods for solving general HJ equations may not exist. 
We give two preliminary definitions.\\

\begin{dfn}[Safe node] \label{def:safe}
Let $x_i\in$ \CONSsp and let $x_{i,1}^*,...,x_{i,p}^*$ be the neighboring interpolation points of $x_i$ 
achieving the minimum in \eqref{SLscheme} ($p=2$ or $p=3$ depending on the employed SL scheme). Denote by $\lambda_{i,1}^*,...,\lambda_{i,p}^*$ the corresponding interpolation weights and define, for $j=1,...,p$, 
$$
b_{i,j}=
\left\{
\begin{array}{ll}
1 & \mbox{if } x^*_{i,j}\in \mbox{\ACC} \\
0 & \mbox{otherwise.}
\end{array}
\right.
$$
The node $x_i$ is said to be \emph{safe} if $\sum_{j=1}^p \lambda^*_{i,j} b_{i,j}=1$.\\
\end{dfn}

\begin{figure}[h!]
\begin{center}
\begin{tabular}{ccc}
\includegraphics[width=.35\textwidth]{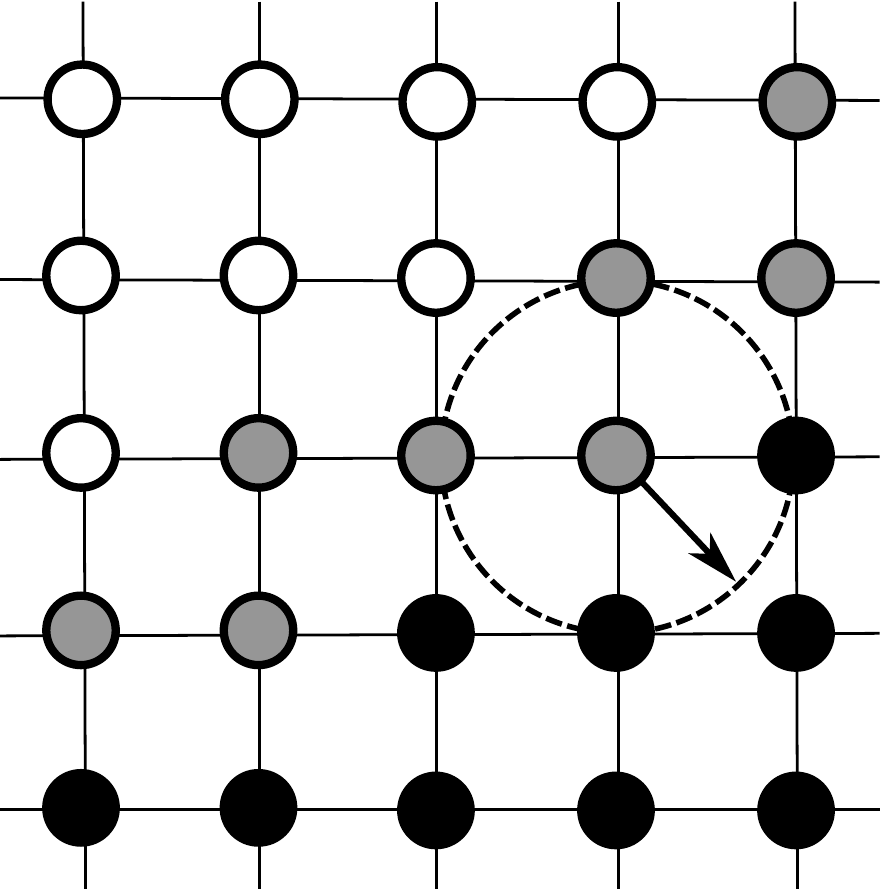}
& \qquad \qquad &
\includegraphics[width=.35\textwidth]{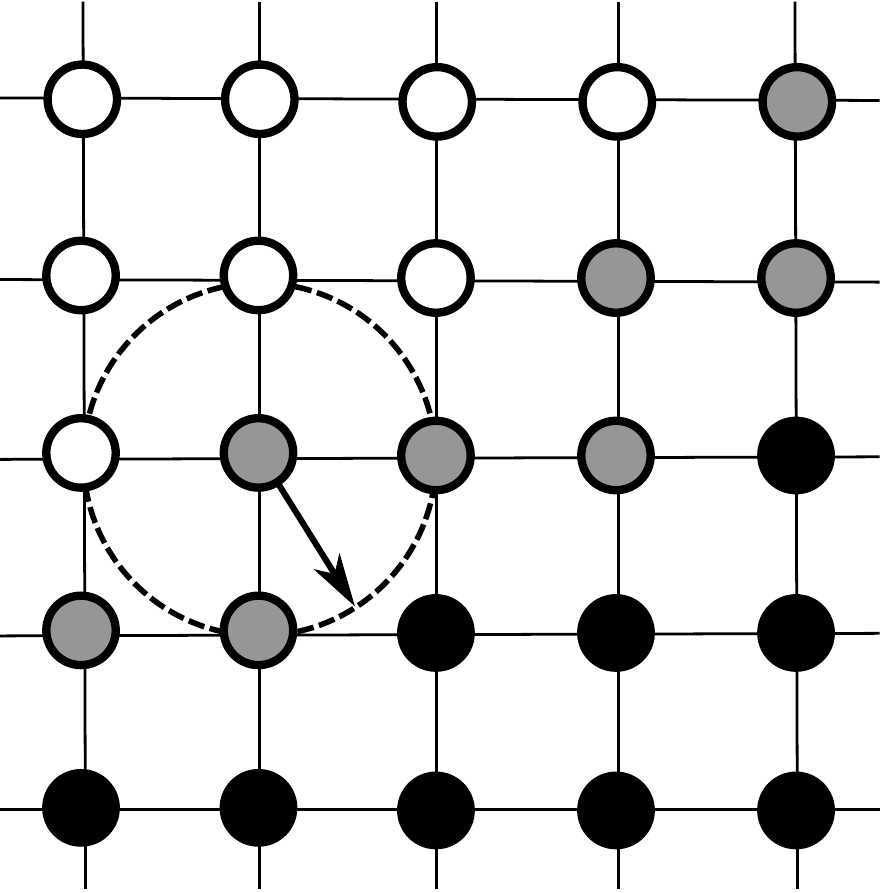}\\
(a) & \qquad \qquad &  (b)
\end{tabular}
\put(-230,0){\small $x_i$}
\put(-234,-27){\small $x_{i,1}^*$}
\put(-188,-27){\small $x_{i,2}^*$}
\put(-188,0){\small $x_{i,3}^*$}
\put(-108,0){\small $x_i$}
\put(-113,-27){\small $x_{i,1}^*$}
\put(-68,-27){\small $x_{i,2}^*$}
\put(-68,0){\small $x_{i,3}^*$}
\end{center}
\caption{Black nodes belong to \ACC, grey nodes belong to \CONS, and white nodes belong to \FAR. 
The arrow denotes the optimal velocity field $f(x_i,a^*)$, and  $x_{i,1}^*,x_{i,2}^*,x_{i,3}^*$ 
are the corresponding neighboring interpolation points involved in the computation of $T(x_i)$.
(a) $x_i$ is safe: its value depends on nodes in \ACCsp only. (b) $x_i$ is not safe: its value depends on nodes in \ACCsp and \CONS. }
\label{fig:safe}
\end{figure}
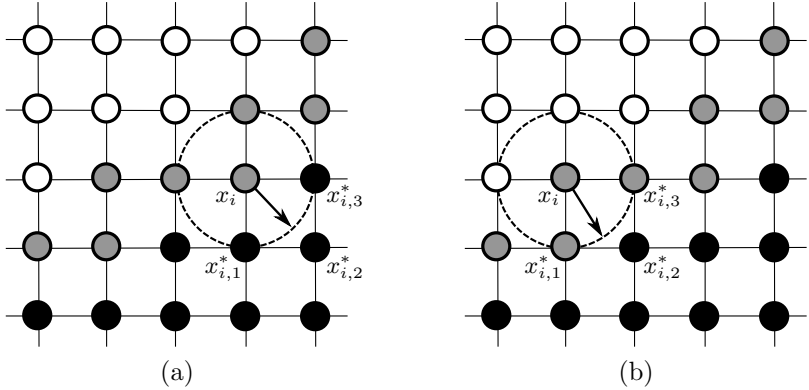

\noindent The previous definition (see Fig.\ \ref{fig:safe} in the case $p=3$), 
means that the computation at $x_i$ involves values of nodes in \ACCsp only. 
By Remark \ref{rem:includeACC}, it is clear that the notion of \emph{safeness} makes sense only if the scheme in use can employ nodes in \CONS. 
Otherwise, all the nodes are safe by construction. Then, we always allow the scheme to use nodes in \CONS.
From a practical point of view, we remark that round-off errors can prevent the safeness to be satisfied. Then, we can relax this condition by requiring that $\sum_{j=1}^p \lambda_{i,j} b_{i,j}\ge1-\sigma$ for some tolerance $\sigma>0$.\\

\begin{dfn}[Exact node] \label{def:exact}
Let $T^{exact}$ be the solution computed by ITM setting $tol$ equal to machine precision in \eqref{stopcriterion}. A node $x_i\in$ \CONSsp is said to be \emph{exact} if its value $T(x_i)$ coincides with $T^{exact}(x_i)$ up to the machine precision.\\
\end{dfn}

\noindent The two new FM-like tools are:\\

\paragraph{Safe Fast Marching Method (SFMM)} This method is identical to FMM but for the rule of acceptance: at each step, the safe node with minimal value in \CONSsp enters \ACC. This method is local and single-pass. 

\medskip

\paragraph{Safe Method (SM)} This method is identical to FMM but for the rule of acceptance: at each step, whichever safe node in \CONSsp enters \ACC. This method is local and single-pass.

Note that the existence of safe nodes in \CONSsp is not in general guaranteed, then SFMM and SM could stop prematurely before \ACCsp covers the whole domain. 
In the next section we will discuss when these method can be successfully employed. 

\medskip
\noindent The two verification methods are:\\

\paragraph{Safe Dumb Method (SDM)} This method is identical to FMM but for the rule of acceptance: at each step, whichever safe and exact node in \CONSsp enters \ACC. 

This method is not usable in practice, since it assumes that one already knows the solution of ITM (or any other equivalent method). 
It clearly computes the same solution of ITM and it is nonlocal, because it employs the information contained in the final solution of ITM, defined everywhere. 
SDM is introduced for theoretical purposes since it represents a limit for the applicability of any local single-pass method. 
Indeed, if there is no nodes in \CONSsp which are both safe and exact, we can conclude that the numerical domain of dependence of every exact node in \CONSsp includes nodes in \CONS. Since one cannot say if values at nodes in \CONSsp are exact or not, we face a loop dependency that cannot be resolved keeping the method local and single-pass. As a consequence, if SDM fails, then any local single-pass method will fail.

\medskip
\paragraph{Dumb Method (DM)} This method is identical to FMM but for the rule of acceptance: at each step, whichever exact node in \CONSsp enters \ACC. Similar considerations made for SDM apply. 

This method is also introduced for theoretical purposes. 
Unexpectedly, DM does not always work. For some pathological dynamics and choice of the mesh, it can happen that \CONSsp does not contain any exact node, and the algorithm stops (see Section \ref{sec:numtests} for an example). In such cases it seems that enlarging \CONSsp or breaking the single-pass property is the only way to make the algorithm process all the nodes.
\\

\begin{rmk}
Every FM-like method requires a selection rule to move nodes from \CONSsp to \ACC. It is possible that more than one node in \CONSsp satisfies that rule at the same step of the algorithm. In this case we can either move one node at random among the correct ones, or move all the nodes at once. 
In the latter case we often get an additional speed up of the algorithm. For example, in FMM one can find two or more nodes in \CONSsp with the minimum value, 
while in SM one can find two or more safe nodes. Investigating the difference of the two implementations is beyond the scope of the paper, since we are mainly 
interested in the applicability of the methods rather than in their performance. Then, we always move in ACC one node at a time.
\end{rmk}


\section{Applicability of local single-pass methods}
\label{sec:core}
In this section we address the problem of extending the range of applicability of local single-pass methods to general HJ equations. To this end, we focus on three algorithms discussed in the previous sections which are local and single-pass, namely FMM, SFMM and SM.  
In order to point out their features and limitations, we will also employ the two verification methods SDM and DM.

From the numerical point of view, it is meaningful to divide HJB equations into two classes. Given a mesh, we have:\\
\begin{itemize}
\item[(ISO)~] Equations whose characteristic lines coincide or lie in the same simplex of the gradient lines of their solutions. The prototype equation is the eikonal isotropic \eqref{eiko_nonomog}.\\

\item[($\neg$ISO)] Equations for which there exists at least a grid node where the characteristic line and the gradient of the solution do not lie in the same simplex.\\
\end{itemize}

\noindent FMM works for equations of type ISO and fails for equations of type $\neg$ISO (see \cite{SV03} for further details and explanations).
Let us introduce two other classes for HJB equations of type (\ref{sect1:HJB}):\\

\begin{itemize}
\item[(REG)~] Equations with non-crossing (regular) characteristic lines. Characteristics spread from the target $\mathcal T$ to the rest of the domain without intersecting.\\

\item[($\neg$REG)] Equations with crossing characteristic lines. Characteristics start from the target $\mathcal T$ and then meet in finite time, creating shocks. As a result, the solution $T$ is not differentiable at shocks.
\end{itemize}

\medskip

\noindent Let us comment the applicability of the local single-pass methods by making use of the classifications introduced above.

\medskip \paragraph{(1) SM solves REG}\label{C1}
SM can be applied in the case REG, provided SDM works. Let us denote by $x_i$ one safe node in \CONSsp ($x_i$ exists because we assume SDM can be applied). 
By definition of safeness, the value at $x_i$ only depends on values at nodes in \ACC, that can be assumed to be exact by induction. 
Then, the exactness of the value at $x_i$ is guaranteed by the property REG, which implies that no characteristics will reach $x_i$ 
in the future from another direction, possibly changing its value. 
In other words, the information passes through $x_i$ one and only one time. Then, once $x_i$ is reached by the region \ACC, it is ready to enter \ACC.

\medskip
\paragraph{(2) Is the minimum-value rule really needed?\hspace{-4.15pt}}
Having in mind the FMM (and its ancestor, the Dijkstra's algorithm \cite{D59}), one can be convinced that giving priority to the smallest value among nodes in \CONSsp is an essential request to 
make the method work. On the contrary, by the above comment \textit{(1)}, 
we know that a method like SM, which makes no distinction among nodes with respect to their values, works in the case REG (both ISO and $\neg$ISO), provided SDM works. 
The choice of the minimum value becomes essential only in the $\neg$REG case, where characteristics reach some point from two or more different directions. 
We discuss this point in the next comment. 

\medskip
\paragraph{(3) Handling shocks in the $\neg$REG case} 
Let us consider the $\neg$REG case and let $x$ be a point belonging to a shock, i.e.\ where the solution is not differentiable. By definition, the value $T(x)$ is carried by two or more characteristic lines reaching $x$ at the same time. Similarly, let $x_i$ be a grid node $\Delta x$-close to the shock. In order to mimic the continuous case,
$x_i$ has to be approached by the \ACCsp region approximately at the same time from the directions corresponding to the characteristic lines. 
In this case, the value $T(x_i)$ is correct (no matter the upwind direction is chosen) and, more important, the characteristic information stops at $x_i$ and it is no longer propagated, getting stuck by the \ACCsp region. 
As a consequence, the shock is localized properly.

We remark that FMM is able to handle shocks when applied to ISO \& $\neg$REG equations. Indeed, thanks to the minimum-value rule, the evolving region \CONSsp is, at every time, a good approximation of the level sets of the final solution and shocks are reached by the \ACCsp region approximately at the same step of the algorithm.

\medskip
\paragraph{(4) $\neg$ISO case requires CONS not to be an approximation of the level sets of the solution}
In order to solve correctly $\neg$ISO equations, \CONSsp cannot be at any time an approximation of the level sets of the solution. 
This is due to the fact that the anisotropy shifts the characteristic directions, so that they no more coincide with the gradient directions. 
\CONS, to be correctly enlarged, should not follow the gradient direction and then it coincides no more with the level sets. In Section \ref{sec:numtests} we show an example for equation (\ref{eiko_aniso_omog}). See also \cite{SV03} for a more detailed explanation.

\medskip
\paragraph{(5) Can local single-pass methods solve general HJB equations?\hspace{-4.15pt}}
Let us consider the $\neg$ISO \& $\neg$REG case. By comments \textit{(3)}-\textit{(4)}, in order to solve $\neg$ISO equations, the \CONSsp region cannot be an approximation 
of a level set of the solution. But, doing so, a node $x_i$ close to a shock can be reached by \ACCsp at different times. 
When \ACCsp reaches $x_i$ for the first time, it is impossible to detect the presence of the shock by using only local information. 
Indeed, only a global view of the solution allows one to know that another characteristic line will reach $x_i$ at a later time. 
As a consequence, the algorithm continues the enlargement of \CONSsp and \ACC, thus making an error that cannot be redressed in the future.
Test 4 in Section \ref{sec:numtests} shows an example in which a shock crosses a region with strong anisotropy. In this situation, it seems impossible to get the correct solution without the addition of nonlocal information regarding the location of the shock, or going back to nodes in \ACCsp at later time. \\




Table \ref{tab:birdeye} summarizes the comments above. Note that the word ``no'' in the table should be meant as ``not in general'', 
since some exceptions are possible. It is plain that SFMM is the more versatile among all methods (whenever it can be applied), since it joins advantages of both SM and FMM.
\begin{table}[h!]
\caption{A bird's eye view on the applicability of local single-pass methods}
\label{tab:birdeye}
\begin{center}
\begin{tabular}{|l|c|c|c|c|}
\hline $ $ & ISO \& REG  &  ISO \& $\neg$REG  &  $\neg$ISO \& REG  &  $\neg$ISO \& $\neg$REG
\\ \hline
\hline FMM  & yes & yes & no                 & no \\
\hline SM   & yes & no  & yes (if SDM works) & no \\
\hline SFMM & yes & yes & yes (if SDM works) & no \\
\hline
\end{tabular}
\end{center}
\end{table}

\section{Numerical tests}\label{sec:numtests}
The first aim of this section is comparing the two semi-Lagrangian schemes described in Section \ref{numschemes}, in order to understand which one has to be preferred for practical implementations of the fast methods described above. 
In the following we will denote the two schemes by SL-2p and SL-3p, respectively. 
The second aim is confirming the theoretical observations in Section \ref{sec:core}, and investigating the theoretical bounds given by SDM and DM, 
in order to understand how much SFMM is close to that limit. 
This will give an idea about how much room is still present for further improvements in the field of local single-pass methods. In all the tests the solution computed by ITM on a  $801^2$ grid (with $tol=10^{-16}$) will be referred to as the ``exact'' solution and will be denoted by $T^{\tiny exact}$. 

In Table \ref{tab:classes} we list five reference HJB equations, together with the class they belong to. In all the cases we set $d=2$, $a=(a_1,a_2)\in B(0,1)$, and $\mathcal T=\{(0,0)\}$. 
Moreover, $\lambda$, $\mu$ and $\varepsilon$ denote generic positive parameters. 
Finally we define $\displaystyle m_{\lambda,\mu}(a)=(1+(\lambda\,a_1+\mu\,a_2)^2)^{-\frac{1}{2}}$ and we denote by $\chi_S$ the characteristic function of a set $S$. 

\begin{table}[h!]
\caption{Equations considered for numerical tests and the class they belong to}
\label{tab:classes}
\begin{center}
\begin{tabular}{|l|l|l|}
\hline Equation & Dynamics & Class
\\ \hline
\hline HJB-A & $f(x,y,a)=a$ & ISO \& REG \\
\hline HJB-B & $f(x,y,a)=(1+\chi_{\{x>1\}})\,a$ & ISO \& $\neg$REG \\
\hline HJB-C & $f(x,y,a)=m_{\lambda,\mu}(a)\,a$ & $\neg$ISO \& REG \\
\hline HJB-D & $f(x,y,a)=(m_{\lambda,\mu}(a)+\varepsilon(x-1)\chi_{\{x>1\}})\,a$ & $\neg$ISO \& $\neg$REG \\
\hline HJB-E & $f(x,y,a)=(1+|x+y|)m_{\lambda,\mu}(a)\,a$ & $\neg$ISO \& $\neg$REG \\
\hline
\end{tabular}
\end{center}
\end{table}

{\bf Test 0 (SL-2p vs\ SL-3p).} In this test we compare the schemes described in Section \ref{numschemes} by means of FSM, 
in terms of accuracy and number of iterations. We consider equations HJB-A and HJB-D (for $\varepsilon=0.02$). 
Relative errors in norm $L^1$ and $L^\infty$ with respect to the ``exact'' solution $T^{\tiny exact}$ are defined as
$$
E_1:=\frac{1}{N}\sum_{i=1}^N \frac{|T^{\tiny exact}(x_i)-\approxT(x_i)|}{|T^{\tiny exact}(x_i)|} \qquad\mbox{and}
\qquad E_\infty:=\max_{i=1,...,N}\frac{|T^{\tiny exact}(x_i)-\approxT(x_i)|}{|T^{\tiny exact}(x_i)|}\,.
$$
By ``sweep'' we mean \emph{four} iterations executed in four different directions. When reporting the number of sweeps of FSM, we include the final ``stopping'' sweep, needed to realize that convergence is reached, namely the stopping rule (\ref{stopcriterion}) is satisfied. We choose $tol=10^{-16}$ (machine precision). 
Results are reported in Table \ref{tab:SL2p_vs_SL3p}.
\begin{table}[h!]
\caption{Test 0: SL-2p and SL-3p schemes comparison}
\label{tab:SL2p_vs_SL3p}
\begin{center}
\begin{tabular}{|l|c|c|c|c|c|c|}
\hline equation & grid & scheme & $E_\infty$ & $E_1$  &  \# sweeps
\\ \hline
\hline HJB-A & $101^2$ & SL-2p & 0.130 & 0.016 & $3$ \\
\hline HJB-A & $101^2$ & SL-3p & 0.079 & 0.009 & $3$ \\
\hline HJB-A & $201^2$ & SL-2p & 0.094 & 0.008 & $3$ \\
\hline HJB-A & $201^2$ & SL-3p & 0.058 & 0.004 & $3$ \\
\hline HJB-A & $401^2$ & SL-2p & 0.050 & 0.003 & $3$ \\
\hline HJB-A & $401^2$ & SL-3p & 0.030 & 0.002 & $3$ \\
\hline
\hline HJB-D & $101^2$ & SL-2p & 0.888 & 0.053 & $8$ \\
\hline HJB-D & $101^2$ & SL-3p & 0.635 & 0.029 & $4$ \\
\hline HJB-D & $201^2$ & SL-2p & 0.535 & 0.027 & $7$ \\
\hline HJB-D & $201^2$ & SL-3p & 0.405 & 0.014 & $4$ \\
\hline HJB-D & $401^2$ & SL-2p & 0.245 & 0.010 & $7$ \\
\hline HJB-D & $401^2$ & SL-3p & 0.189 & 0.005 & $3$ \\
\hline
\end{tabular}
\end{center}
\end{table}

We recall that, as discussed in Section \ref{algorithms}, the convergence of FSM is ensured in 1 sweep for equation HJB-A, see \cite{Z05}. Nevertheless, real algorithms involving double precision computations can require 2 sweeps to reach the machine precision. 
The third sweep reported in Table \ref{tab:SL2p_vs_SL3p} is the ``stopping'' sweep. \\
It is rather clear that SL-3p overcomes SL-2p in terms of both accuracy and number of sweeps. 
This is likely due to the fact that SL-3p can propagate the characteristic information of the HJB equation along diagonal directions easier than SL-2p. 

Dealing instead with FM-like methods, the two schemes show a difference in the order of acceptance of the nodes in \CONS. In particular, we noted that algorithms based on SL-3p provide a larger number of safe nodes in \CONS, thus extending the applicability of SM, SFMM, and SDM.

From now on, only the scheme SL-3p will be employed for all the following tests.

\medskip

{\bf Test 1 (ISO \& REG).} In this test we compare FSM, FMM and SM against HJB-A. 
Errors with respect to the ``exact'' solution  $T^{\tiny exact}$ are reported in Table \ref{tab:eik_diff}.
\begin{table}[h!]
\caption{Test 1: ISO \& REG}
\label{tab:eik_diff}
\begin{center}
\begin{tabular}{|c|l|c|c|c|}
\hline grid & method & $E_\infty$ & $E_1$
\\ \hline
\hline $101^2$ & FSM & 0.079 & 0.009 \\
\hline $101^2$ & FMM & 0.079 & 0.009 \\
\hline $101^2$ & SM  & 0.079 & 0.009 \\
\hline
\hline $201^2$ & FSM & 0.057 & 0.004 \\
\hline $201^2$ & FMM & 0.057 & 0.004 \\
\hline $201^2$ & SM  & 0.057 & 0.004 \\
\hline
\hline $401^2$ & FSM & 0.029 & 0.001 \\
\hline $401^2$ & FMM & 0.029 & 0.001 \\
\hline $401^2$ & SM  & 0.029 & 0.001 \\
\hline
\end{tabular}
\end{center}
\end{table}\\
The three methods lead to the same error because they compute exactly the same solution. 
A fortiori, also SFMM does. This confirms that SM can be applied in ISO \& REG cases and that picking the minimum value in \CONS, as acceptance rule, 
is not strictly needed here to compute the correct solution. 

\medskip

\begin{figure}[h!]
\begin{center}
\begin{tabular}{cc}
\includegraphics[width=0.4\textwidth]{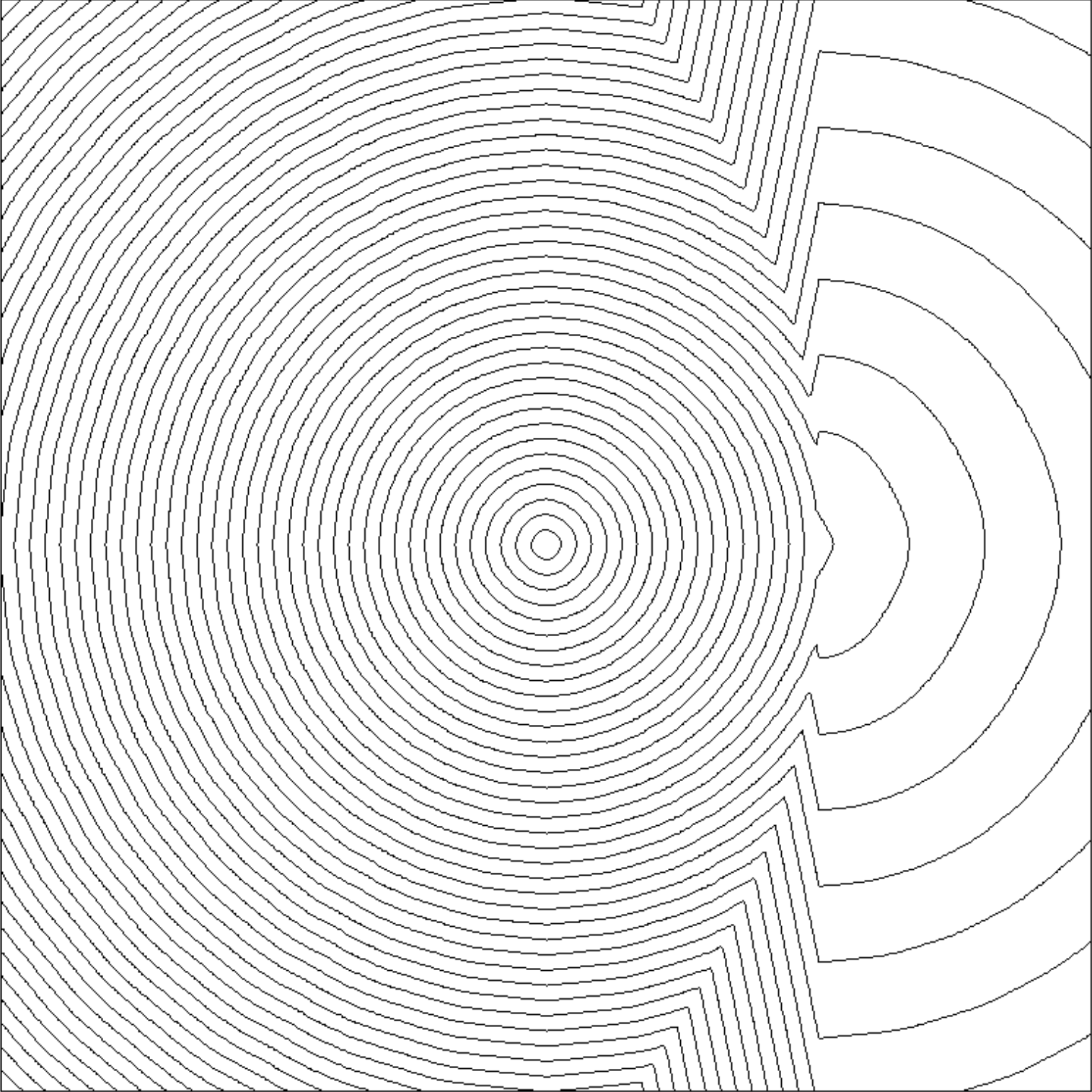} &
\includegraphics[width=0.4\textwidth]{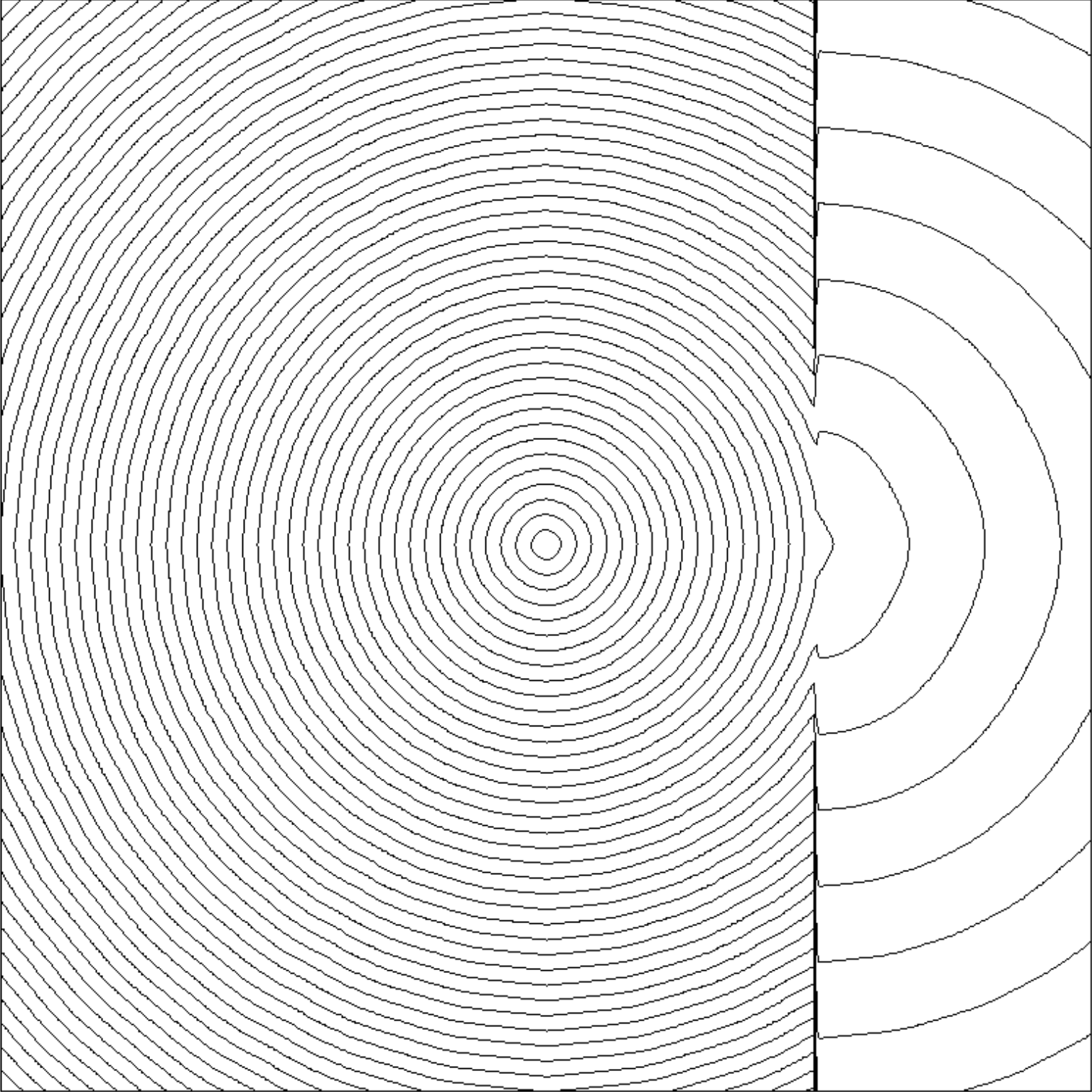} \\
FSM & SM
\end{tabular}
\caption{Test 2: level sets of the solutions computed by FSM and SM}
\label{fig:eik_!diff}
\end{center}
\end{figure}

{\bf Test 2 (ISO \& $\neg$REG).} In this test we compare FSM, FMM, SFMM and SM against HJB-B. Errors with respect to the ``exact'' solution  $T^{\tiny exact}$ are reported in 
Table \ref{tab:eik_!diff}. 
\begin{table}[h!]
\caption{Test 2: ISO \& $\neg$REG}
\label{tab:eik_!diff}
\begin{center}
\begin{tabular}{|c|l|c|c|c|}
\hline grid & method & $E_\infty$ & $E_1$
\\ \hline
\hline $101^2$ & FSM  & 0.079 & 0.011 \\
\hline $101^2$ & FMM  & 0.079 & 0.011 \\
\hline $101^2$ & SFMM & 0.079 & 0.011 \\
\hline $101^2$ & SM   & 0.583 & 0.019 \\
\hline
\hline $201^2$ & FSM  & 0.057 & 0.006 \\
\hline $201^2$ & FMM  & 0.057 & 0.006 \\
\hline $201^2$ & SFMM & 0.057 & 0.006 \\
\hline $201^2$ & SM   & 0.606 & 0.014 \\
\hline
\hline $401^2$ & FSM  & 0.029 & 0.002 \\
\hline $401^2$ & FMM  & 0.029 & 0.002 \\
\hline $401^2$ & SFMM & 0.029 & 0.002 \\
\hline $401^2$ & SM   & 0.603 & 0.011 \\
\hline
\end{tabular}
\end{center}
\end{table}
FSM, FMM and SFMM lead to the same error because they compute exactly the same solution. Conversely, 
SM cannot be used here, since it is not able to properly locate the shocks (see Fig.\ \ref{fig:eik_!diff}).

\medskip

{\bf Test 3 ($\neg$ISO \& REG).} In this test we compare FSM, FMM and SM against HJB-C. Errors with respect to the ``exact'' solution  $T^{\tiny exact}$ are reported in Table \ref{tab:!eik_diff}.
\begin{table}[h!]
\caption{Test 3: $\neg$ISO \& REG}
\label{tab:!eik_diff}
\begin{center}
\begin{tabular}{|c|l|c|c|c|}
\hline grid & method & $E_\infty$ & $E_1$
\\ \hline
\hline $101^2$ & FSM  & 0.635 & 0.029 \\
\hline $101^2$ & FMM  & 0.635 & 0.058 \\
\hline $101^2$ & SM   & 0.635 & 0.029 \\
\hline
\hline $201^2$ & FSM  & 0.404 & 0.014 \\
\hline $201^2$ & FMM  & 0.408 & 0.049 \\
\hline $201^2$ & SM   & 0.404 & 0.014 \\
\hline
\hline $401^2$ & FSM  & 0.189 & 0.005 \\
\hline $401^2$ & FMM  & 0.290 & 0.044 \\
\hline $401^2$ & SM   & 0.189 & 0.005 \\
\hline
\end{tabular}
\end{center}
\end{table}
FSM and SM lead to the same error because they compute exactly the same solution. 
A fortiori, also SFMM does. Conversely, FMM fails (although it is quite robust), since it does not compute the same solution of FSM. 
This comes from the fact that FMM is not able to deal with substantial anisotropies, as discussed in Section \ref{algorithms} (see also \cite{SV03} for more details).

\medskip

\begin{figure}[h!]
\begin{center}
\includegraphics[width=0.4\textwidth]{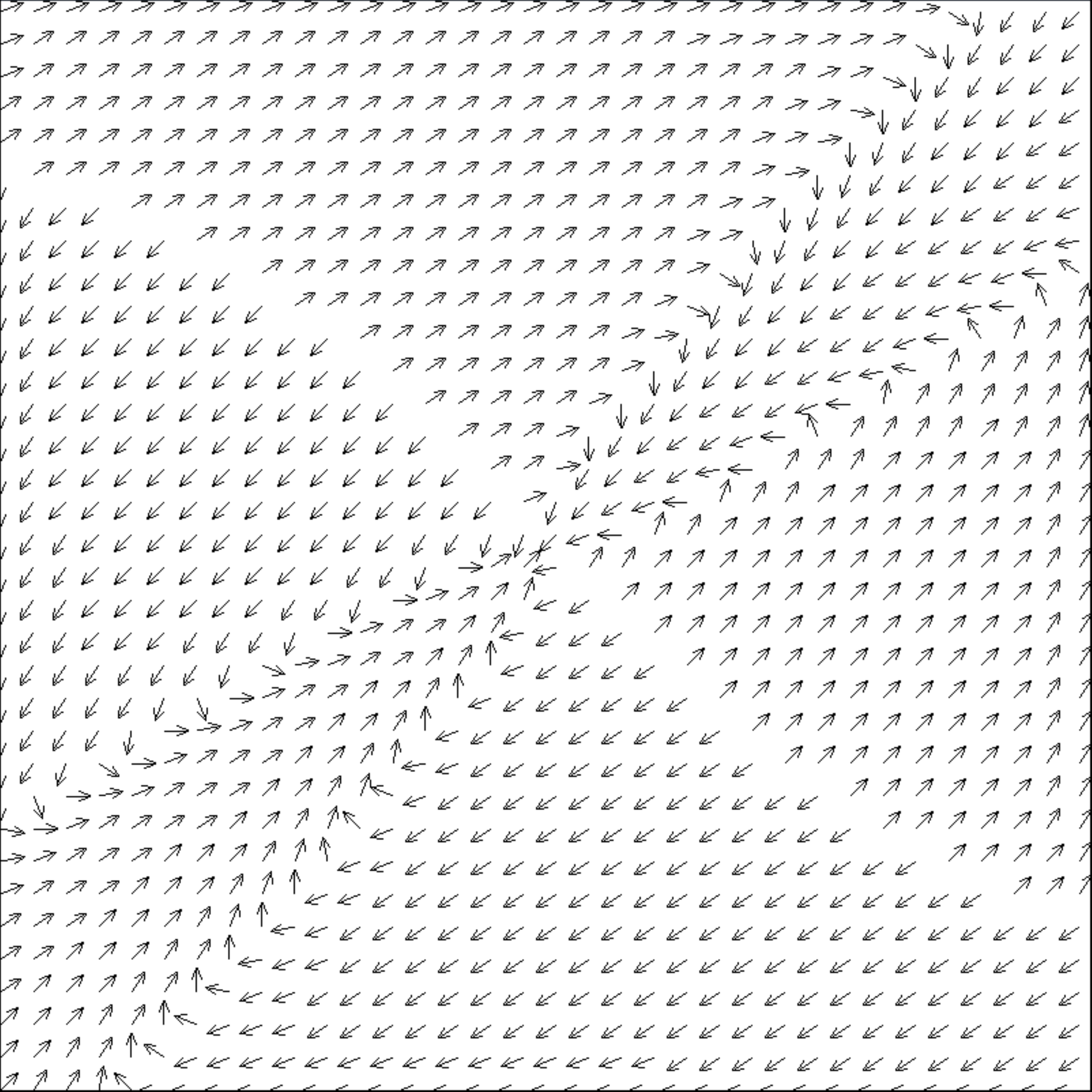}
\caption{Test 4: Optimal vector field $f^*_{\textup{FSM}}$ computed by FSM (zoom around the origin)}
\label{fig:!eik_!diff_frecce}
\end{center}
\end{figure}

\begin{figure}[h!]
\begin{center}
\begin{tabular}{cc}
\includegraphics[width=0.4\textwidth]{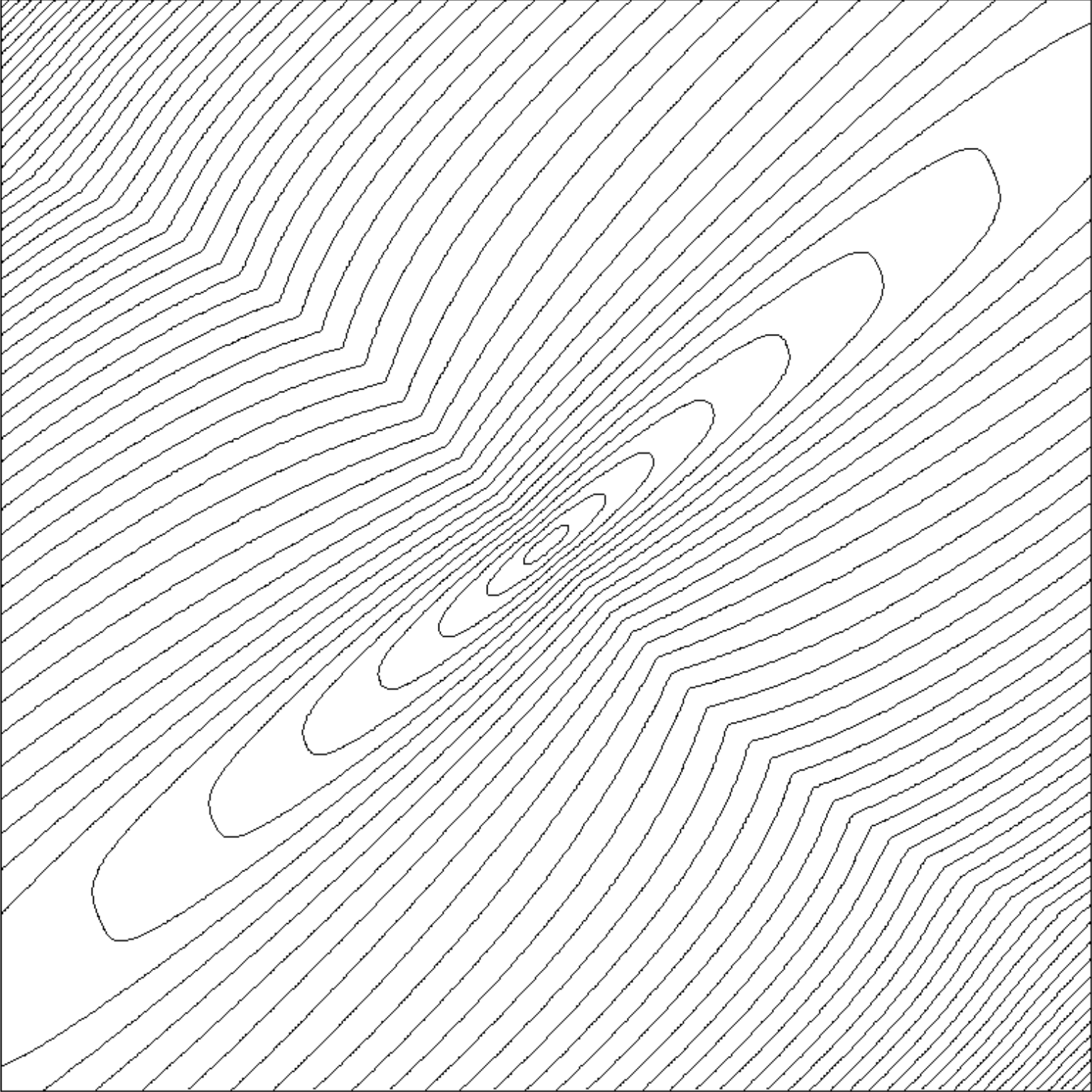} &
\includegraphics[width=0.4\textwidth]{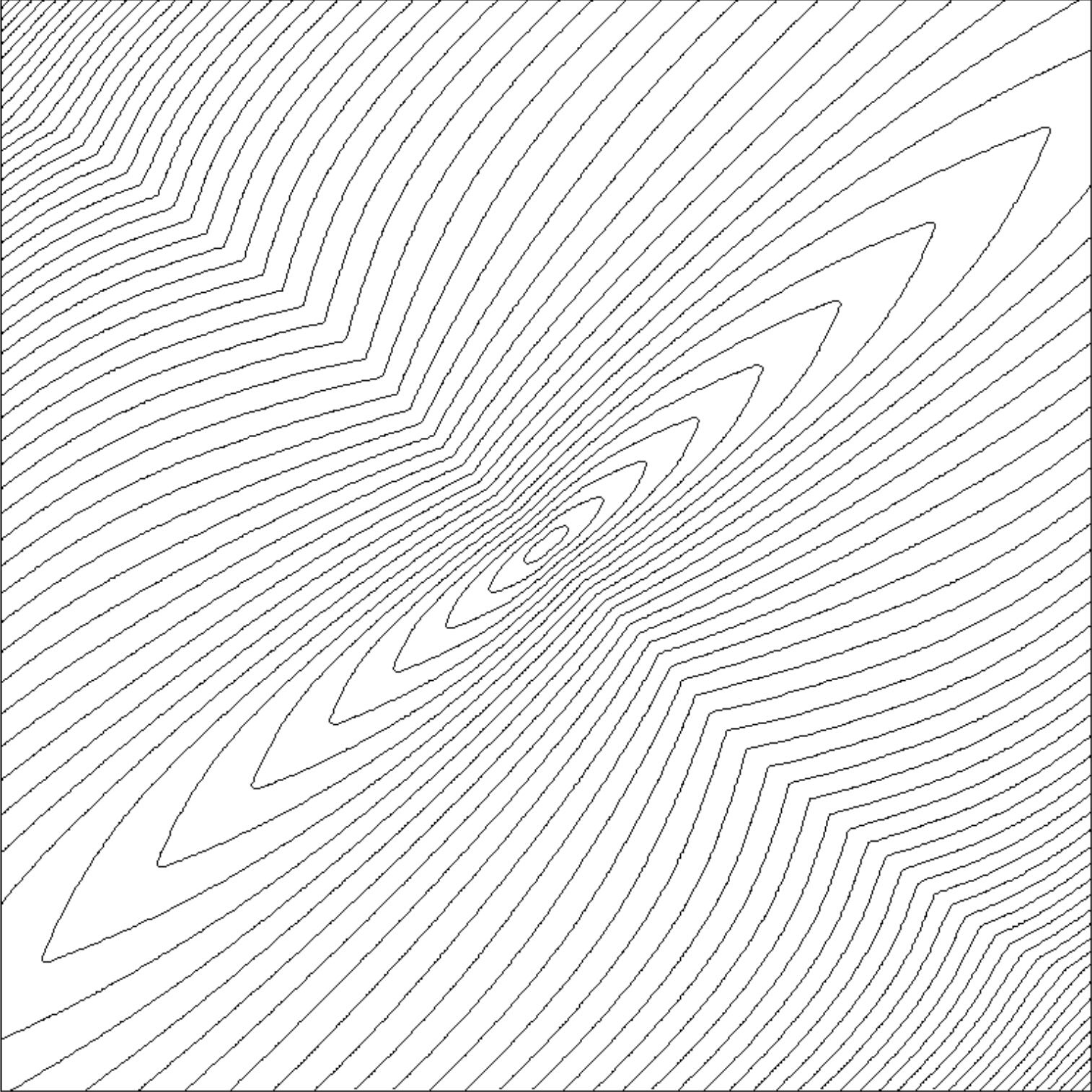}\\
FSM & FMM \\
\includegraphics[width=0.4\textwidth]{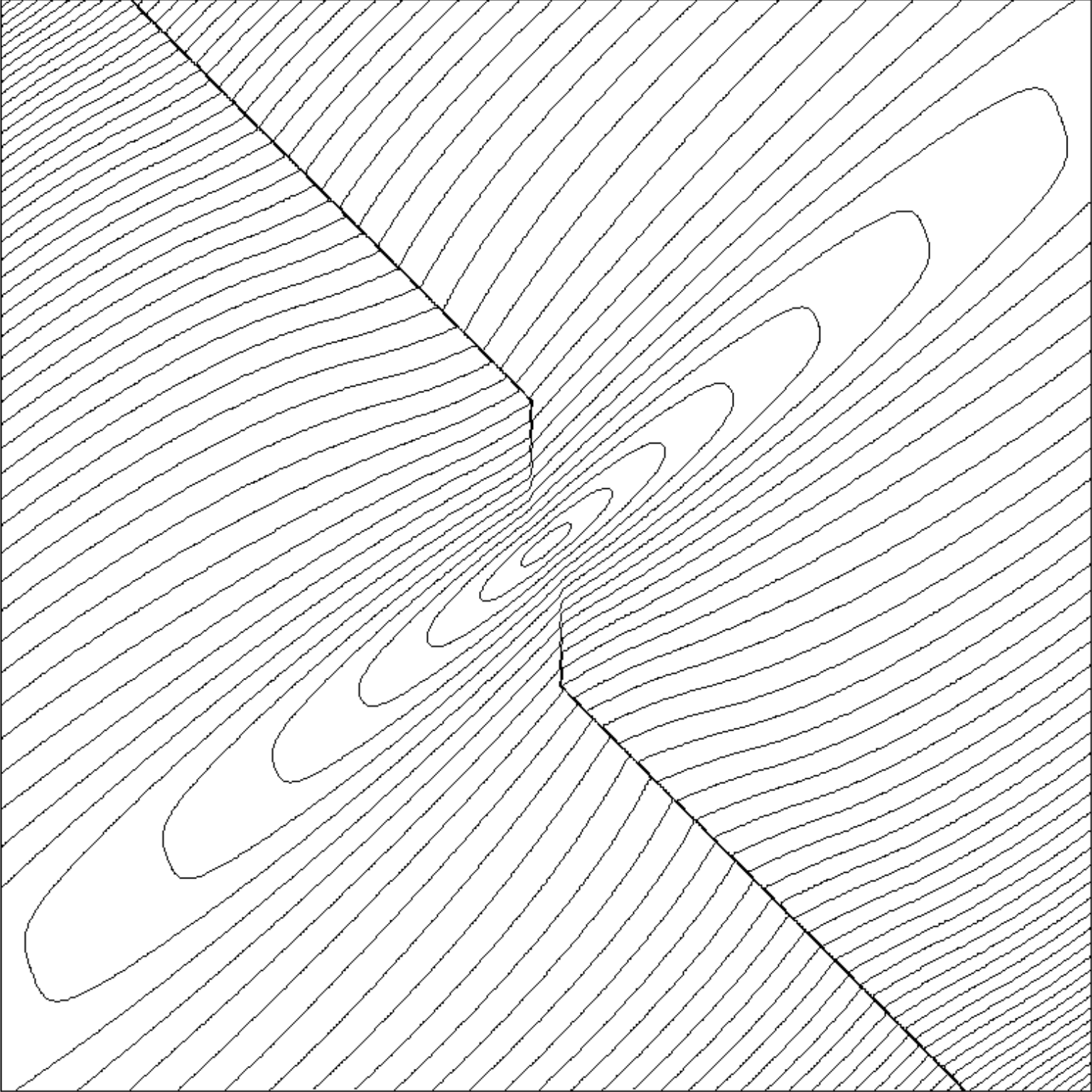} &
\includegraphics[width=0.4\textwidth]{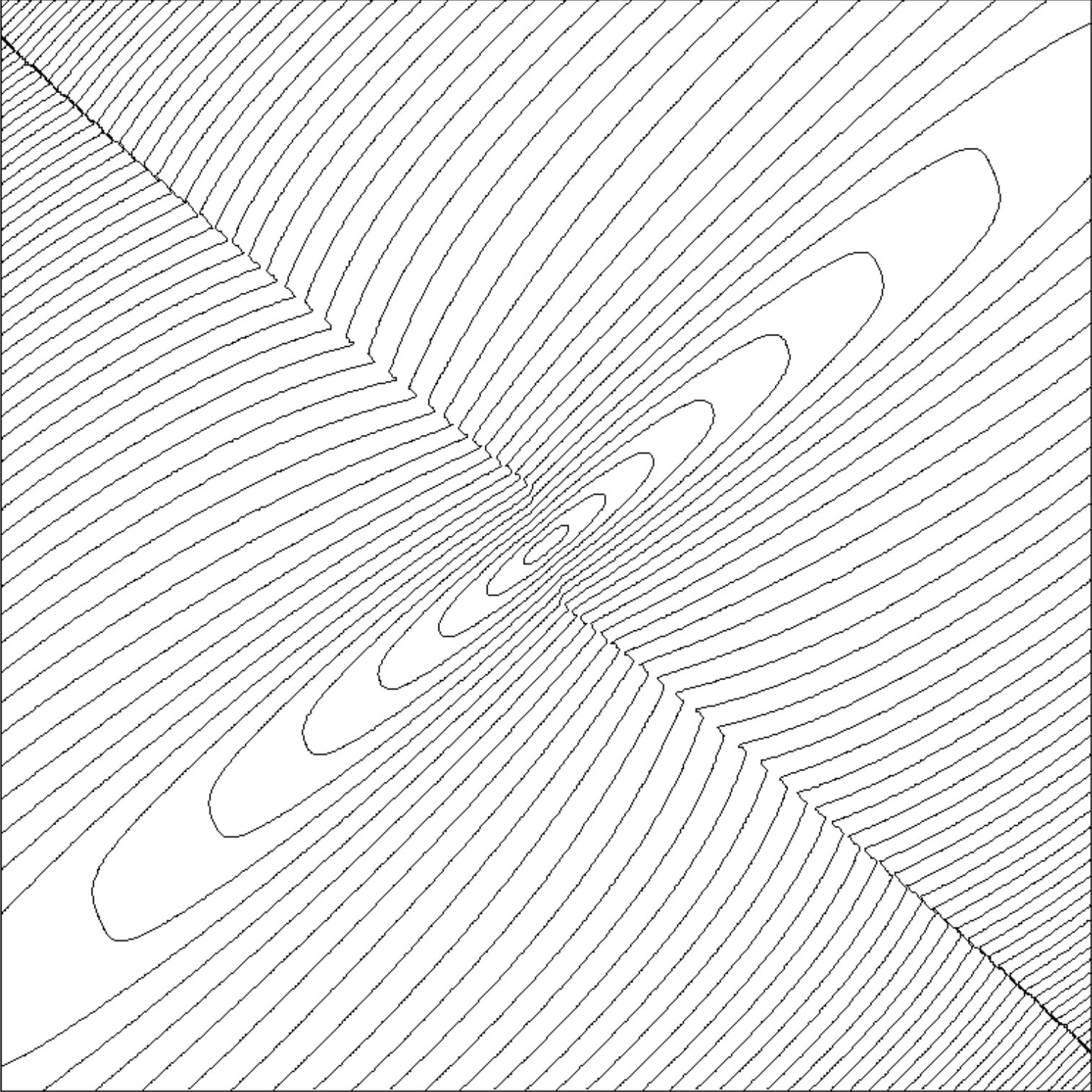} \\
SM & SFMM
\end{tabular}
\caption{Test 4: level sets of the solutions computed by FSM, FMM, SM and SFMM}
\label{fig:!eik_!diff}
\end{center}
\end{figure}

{\bf Test 4 ($\neg$ISO \& $\neg$REG).} In this test we compare FSM, FMM, SFMM, SM, DM and SDM against HJB-E (for $\lambda=6$ and $\mu=5$). 
Fig.\ \ref{fig:!eik_!diff_frecce} shows some optimal directions (characteristic lines) computed by means of FSM.\\\\
Both the strong inhomogeneity (characteristic lines bend hardly in the I and III quadrant) and the shock (the cubic-like curve in the II and IV quadrant) are visible.
\\
Table \ref{tab:!eik_!diff} reports the error $E_1$ with respect to the ``exact'' solution  $T^{\tiny exact}$.
\begin{table}[h!]
\caption{Test 4: $\neg$ISO \& $\neg$REG}
\label{tab:!eik_!diff}
\begin{center}
\begin{tabular}{|c|c|c|c|c|c|c|}
\hline grid/method & FSM & FMM & SFMM & SM & DM & SDM
\\ 
\hline $101^2$ & 0.114 & 0.170 & 0.115 & 0.124 & 0.114 & 0.114 \\
\hline $201^2$ & 0.061 & 0.132 & 0.062 & 0.072 & 0.061 & 0.061 \\
\hline $401^2$ & 0.024 & 0.109 & 0.025 & 0.036 & 0.024 & 0.024\\
\hline
\end{tabular}
\end{center}
\end{table}
In this case only ``dumb'' methods (DM and SDM) are able to compute the same solution of FSM, although SFMM is very close to FSM. 
The differences among the methods are much more evident looking at the level sets of the corresponding solutions, reported in Fig.\ \ref{fig:!eik_!diff}.
FSM is able to respect the anisotropy, indeed the level sets of its solution around the origin are ellipses as expected. Moreover, it properly catches the shock. FMM tries catching the shock, but fails in respecting the anisotropy. 
SM tries respecting anisotropy, but fails in catching the shock. Finally, SFMM is a kind of mix between FMM and SM.

\medskip

\begin{figure}[h!]
\begin{center}
\begin{tabular}{cc}
\includegraphics[width=0.4\textwidth]{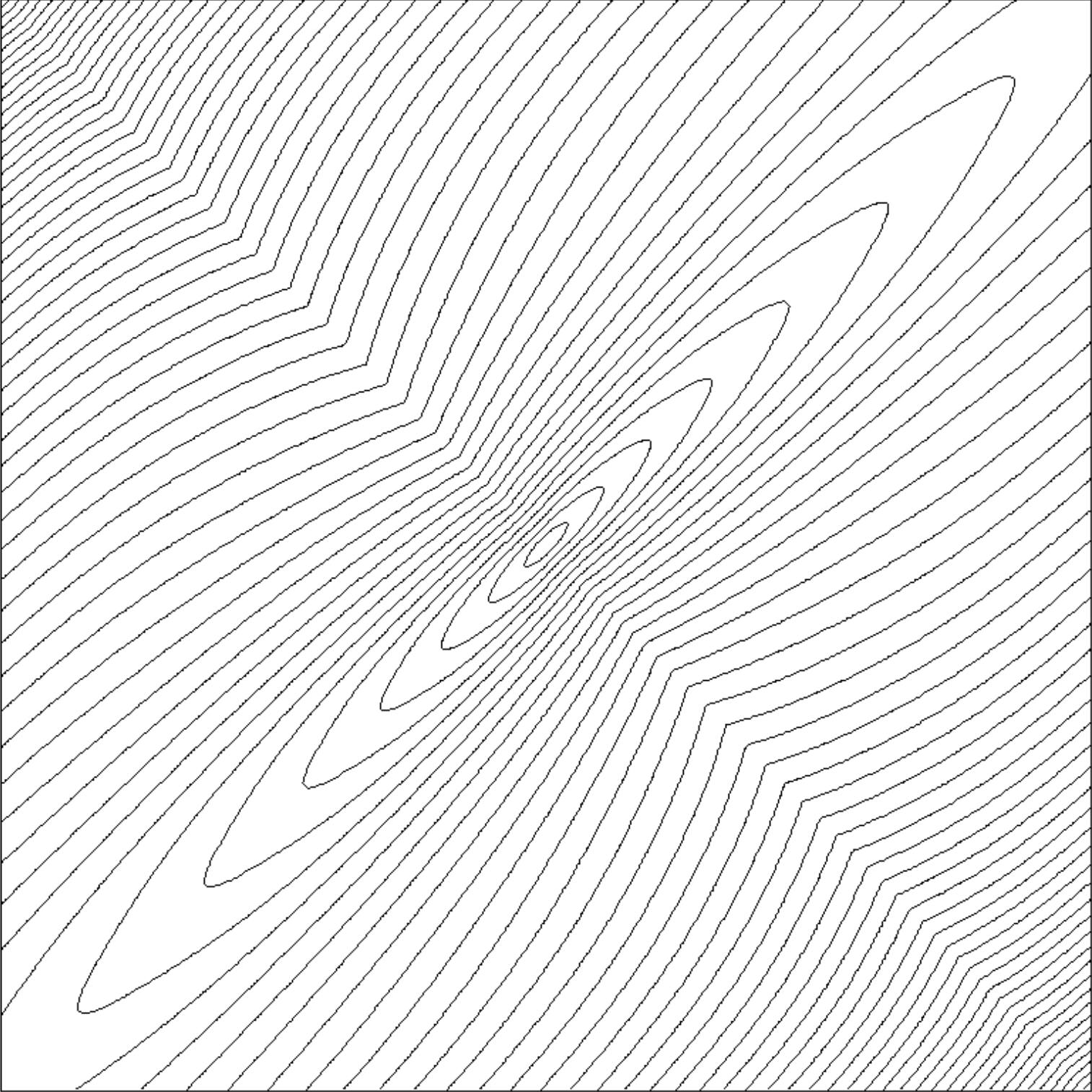}  &
\includegraphics[width=0.4\textwidth]{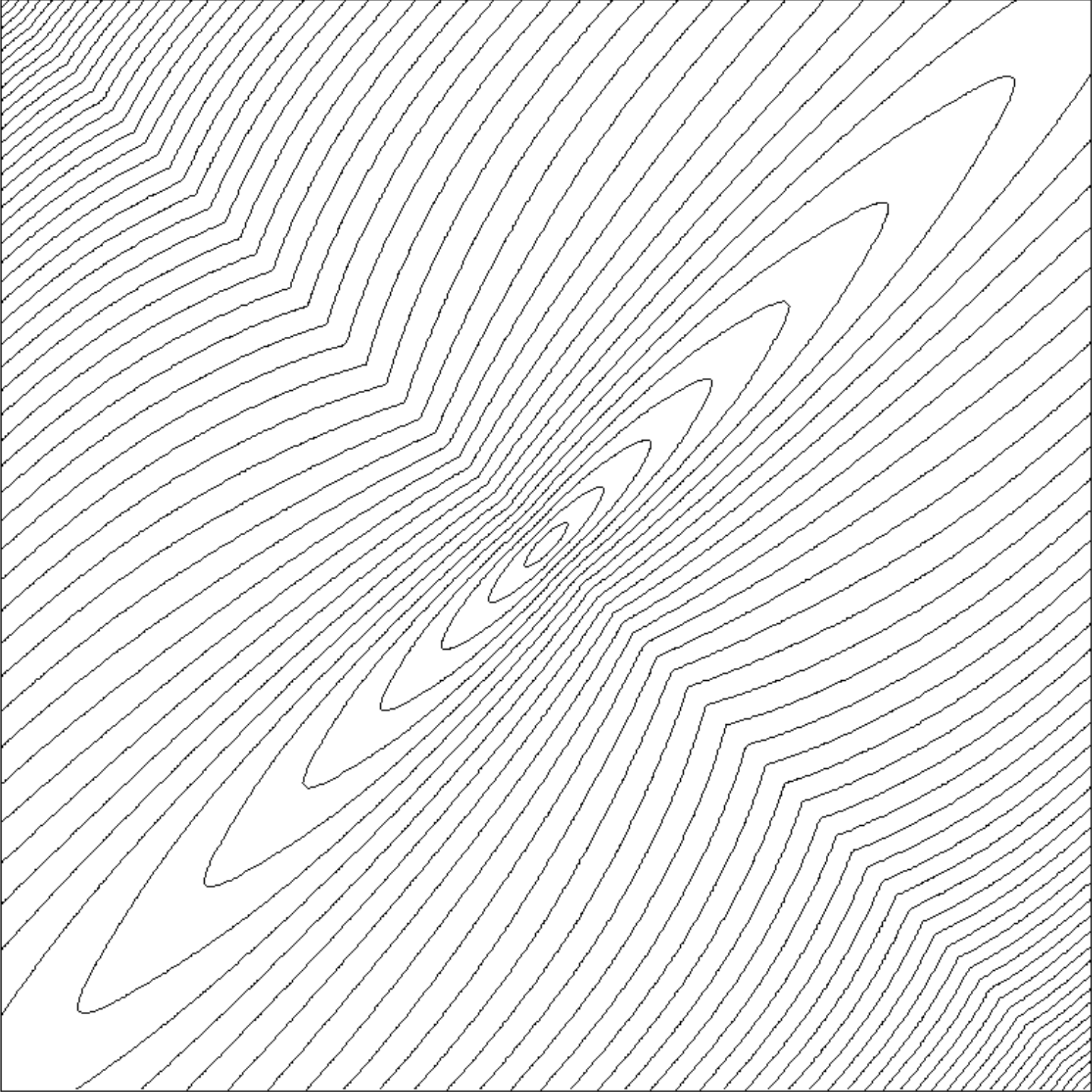} \\
FSM & SFMM 
\end{tabular}
\caption{Test 5: level sets of the solutions computed by FSM and SFMM
}
\label{fig:!eik_!diff_easy}
\end{center}
\end{figure}

{\bf Test 5 ($\neg$ISO \& $\neg$REG - easy case).} In this test  we compare FSM and SFMM against HJB-E (for $\lambda=5$ and $\mu=5$). 
Due to the fact that $\lambda=\mu$, the shock has a particular symmetry with respect to the axes. 
This symmetry makes SFMM work, since \CONSsp ``luckily'' reaches the shock at the same time from both sides (see Fig.\ \ref{fig:!eik_!diff_easy}).
This example shows that local single-pass schemes can solve $\neg$ISO \& $\neg$REG equations in some special cases.

\medskip

\begin{figure}[h!]
\begin{center}
\includegraphics[width=.9\textwidth]{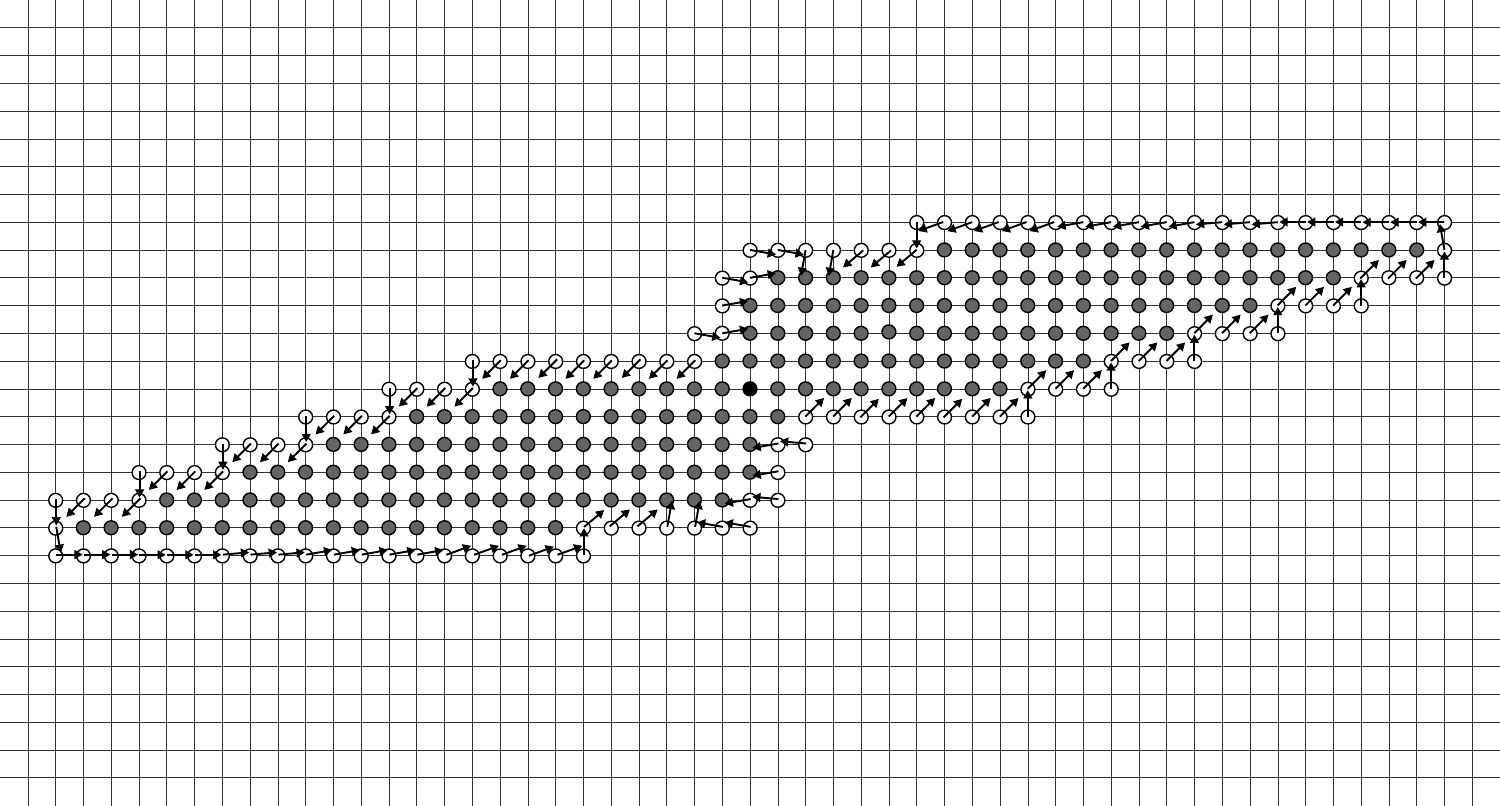}
\caption{Test 6: SDM fails because no safe nodes are found in \CONS. Note the loop dependency following the optimal vector field $f^*_{\textup{SDM}}$}
\label{fig:!eik_!diff_hard}
\end{center}
\end{figure}

{\bf Test 6 ($\neg$ISO \& $\neg$REG - hard case).} In this test we show that even SDM and DM can fail in computing the correct solution, 
i.e.\ it can happen that either there is no safe node in \CONSsp and/or there is no exact node in \CONS. 
Therefore, the methods stop abruptly before \ACCsp covers the whole domain.
We consider again the equation HJB-E (for $\lambda=10$ and $\mu=5$). This case is even more pathological than that depicted in Fig.\ \ref{fig:!eik_!diff_frecce}: 
characteristic lines bend too much compared to the mesh size, i.e. they can significantly change direction within a single cell.
We caught the precise moment in which both SDM and DM stop working, due to the lack of safe and exact nodes in \CONS. 
In Fig.\ \ref{fig:!eik_!diff_hard} the black central node is the target, gray nodes represent the \ACCsp region, whereas white nodes are in \CONS.
For each node in \CONSsp we plot the optimal vector field $f^*_{\textup{SDM}}$ computed by means of the current solution of SDM. 
It is evident that every node in \CONSsp depends on other nodes in \CONS, so that a loop is created and no safe node is present.

\begin{figure}[h!]
\begin{center}
\includegraphics[width=.45\textwidth]{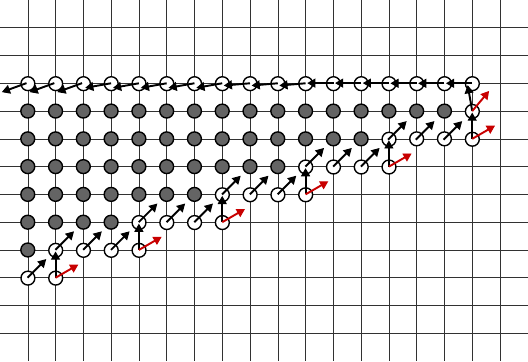}
\caption{Test 6: DM fails because no exact nodes are found in \CONS. Exact information flows back from \FARsp region breaking locality and single-pass property}
\label{fig:!eik_!diff_hard_zoom}
\end{center}
\end{figure}

In Fig.\ \ref{fig:!eik_!diff_hard_zoom} we show a detail of Fig.\ \ref{fig:!eik_!diff_hard} and we plot the optimal vector fields  $f^*_{\textup{DM}}$ (in black) and
 $f^*_{\textup{ITM}}$ (in red), computed by means of the solution of DM and ITM respectively (if the two optimal vector fields coincide only one is plotted). 
\\
It can be seen that $f^*_{\textup{ITM}}$ points either toward the \FARsp region or toward nodes in \CONS, whose $f^*_{\textup{ITM}}$ also points toward the \FARsp region (recursively). 
This means that the values at nodes in \ACCsp and \CONSsp are not enough to compute exact values in \CONS, even if we perform an additional stabilization by iterating the scheme on \CONSsp up to convergence. We infer that we are facing a large loop in the numerical domain of dependence of the nodes in \CONS, which includes also nodes currently in \FAR. This particular behavior of the characteristic flow is also confirmed by the fact that ITM requires in this case a huge number of iterations to reach convergence, compared to that of the previous tests.   
The thickness of \CONSsp must be increased (as in BFMM), in order to resolve the dependency. 
\section{Conclusions}
The above tests and considerations allow us to sketch some final comments and suggest new  directions for a future analysis of acceleration methods. First of all, we want to stress that all the considerations debated in the paper have a theoretical value. Indeed, from the  practical point of view, it is not always possible to know in advance if an equation falls in the class ISO or REG and then it is not evident how to choose a method which is able to solve it. Only ITM and FSM can be safely used if no \textit{a priori} knowledge of the solution is available.

1. SM is one of the simplest methods one can imagine, nevertheless it is able to solve a large class of equations, including the homogeneous anisotropic eikonal equation (\ref{eiko_aniso_omog}). Therefore, methods such OUM and PFMM are in some sense ``too complicated'' than necessary. In our opinion, the reason why the minimum-value rule has been given a crucial role so far is simply that Dijkstra method uses it. Nevertheless, on graphs, the distinction between REG and $\neg$REG is not visible, as well as the condition of safeness. The importance of the latter condition was missed because the different possibilities to use \CONSsp nodes  has been completely underestimated (see Remark \ref{rem:includeACC}).
In this respect, we point out that SM is very similar to CFMM since the acceptance rule used in CFMM actually coincides with that of SM (see Definition \ref{def:safe}). 
In \cite{CF08} it was already noted that, running CFMM, \CONSsp does not coincide with the level set of the solution, but this fact was not fully exploited as we have done in this paper. 

2. The reliability of the SFMM to solve  $\neg$ISO \& $\neg$REG equations cannot be known in advance. Even if  the method computes a solution, i.e.\ \ACCsp covers the whole grid, we do not know if that solution can be correct or not. On the contrary, if the method stops, due to the lack of safe nodes in \CONS,  the user can eventually conclude that this method cannot be used for the equation under consideration. 

3. Our experience suggests that there is no much room between SFMM and SDM,  meaning that it is quite difficult to define precisely a class of equations that can be solved by SDM and not by SFMM. 
Since we have seen that SDM is a sort of a limit of applicability for local single-pass schemes, we conclude that it is relatively fruitless investigating new local single-pass schemes. In order to solve more general equations  one should look for new  methods based on larger and dynamic stencils that will likely produce more complicated implementations.

\section*{Acknowledgements}
Authors are grateful to Alexander Vladimirsky and Michael Breu\ss{} for the interesting and motivating discussions they had with them during the preparation of this paper.


\end{document}